\documentclass[onecolumn]{autart}
\usepackage{graphicx}
\input{amssym}

\renewcommand{\theequation}{\arabic{section}.\arabic{equation}}
\begin{document}
\begin{frontmatter}
\title{Group analysis for \\ generalized reaction-diffusion convection equation} 
\thanks[footnoteinfo]{ Corresponding author: Tel. +9821-73913426.
Fax +9821-77240472.}
\author[]{M. Nadjafikhah\thanksref{footnoteinfo}}\ead{m\_nadjafikhah@iust.ac.ir},
\author[]{S. Dodangeh}\ead{s\_dodangeh@mathdep.iust.ac.ir},
\address{School of Mathematics, Iran University of Science and Technology, Narmak, Tehran 1684613114, Iran.}
\renewcommand{\sectionmark}[1]{}
\begin{abstract}
In this paper we discuss about group classification for non-linear
generalized reaction-diffusion convection equation:
$u_t=(f(x,u)u_x)_x+h(x,u)u_x+k(x,u)$, by using Lie-classical
symmetry method. For this, first we find its symmetry group and
then we find differential invariants for resulted group by using
infinitesimal criterion method and at the end reduce modeling
equations by using resulted invariants. We present application of
this group classification in group classification and obtaining
related similarity solution of KPP equation, too.
\end{abstract}
\end{frontmatter}
\section{Introduction}
This paper devoted to group classification of Generalized
Reaction-Diffusion Convection (G-RDC) equation, by using
Lie-classical method.
\begin{eqnarray}
\Delta\;:\;u_t=(f(x,u)u_x)_x+h(x,u)u_x+k(x,u),\label{eq:11}
\end{eqnarray}
Where $u(x,t)$ is unknown function and $f(x,u),h(x,u),k(x,u)$ are
arbitrary functions. The equation (\ref{eq:11}) generalizes a
number of the well known second-order evolution equations,
describing various process in physics, chemistry and biology.
Symmetry group method plays an important role in the analysis of
differential equations. The history of group classification
methods goes back to  Sophus Lie \cite{[9]}. (See
\cite{[2],[3],[4],[8]}). His work devoted to finding symmetry
groups, differential invariants and linearized or reduced equation
for given model. There are several approach to group
classification of differential equation, we apply infinitesimal
method (See \cite{[2],[3],[10]}) for this. There are another
useful articles and accounts about group classification for
similar equations of (\ref{eq:11}) via other methods, (See
\cite{[11],[12],[13]}). In this paper we generalize RDC equation
to G-RDC equation and applay Lie-classical symmetry method via
applied approach. we hope this work be useful to applied and
theoretical readers.
%
%
\section{Group classification for modeling equation}
Let following one-parameter group
\begin{eqnarray}
\overline{x}=x+\varepsilon\xi(x,t,u)+O(\varepsilon^2),\qquad
\overline{t}=t+\varepsilon\eta(x,t,u)+O(\varepsilon^2),\qquad
\overline{u}=u+\varepsilon\varphi(x,t,u)+O(\varepsilon^2),\label{eq:4}
\end{eqnarray}
be symmetry group of modeling equation $\Delta$. We can obtain
$\xi$, $\eta$ and $\varphi$, by using infinitesimal method.

Consider the vector field
$X:=\xi\partial_x+\eta\partial_t+\varphi\partial_u$ in total space
$M=(x,t,u)$ with $p=2$ and $q=1$.  If this vector field be an
infinitesimal generator of $\Delta$'s symmetry group, then
\begin{eqnarray}
X^{(2)}\Delta=0,\qquad \textrm{whenever}\qquad
\Delta=0.\label{eq:2}
\end{eqnarray}
Where $X^{(2)}$ is second prolong of $X$ and has following form:
\begin{eqnarray}
X^{(2)}=X+\varphi^x\partial_{u_x}+\varphi^t\partial_{u_t}+\varphi^{xx}\partial_{xx}+\varphi^{xt}\partial_{u_{xt}}+\varphi^{tt}\partial_{u_{tt}},\label{eq:
1}
\end{eqnarray}
Where $\varphi^x,\varphi^t,\varphi^{xx},\varphi^{xt}$ and
$\varphi^{tt}$ are respectively:
\begin{eqnarray*}
&&\varphi^x=D_xQ+\xi u_{xx}+\eta u_{xt},
\hspace{1.6cm}\varphi^{t}=D_tQ+\xi u_{xt}+\eta u_{tt},\\
&&\varphi^{xx}=D_{xx}Q+\xi u_{xxx}+\eta u_{xxt}, \qquad
\varphi^{xt}=D_{xt}Q+\xi u_{xxt}+\eta u_{xtt},\qquad
\varphi^{tt}=D_{tt}Q+\xi u_{xtt}+\eta u_{ttt}.
\end{eqnarray*}
Where $D_x$, $D_t$ are total derivative with respect to specified
variables, $D_{xx}=D_xD_x$, $D_{xt}=D_xD_t$ and $D_{tt}=D_tD_t$,
and $Q=\varphi-\xi u_x-\eta u_t$ is the corresponding
characteristic of $X$ (See \cite {[2],[3],[4],[8]}). By using
criterion (\ref{eq: 1}), we find:
\begin{eqnarray}
&&\varphi^t=(f_{xx}u_x+f_{xu}u_x^2+f_xu_{xx}+k_x+h_xu_x)\xi+(f_{xu}u_x+f_{uu}u_x^2+f_uu_{xx}+h_uu_x+k_u)\varphi+\nonumber\\
&&\hspace{0.9cm}+(f_x+2f_uu_x+h)\varphi^x+f\varphi^{xx}.\label{eq:3}
\end{eqnarray}
By substituting $\varphi^x,\varphi^t,\varphi^{xx},\varphi^{xt}$
and $\varphi^{tt}$ in (\ref{eq:3}), we have following results:
\begin{center}\small  \mbox{ } \hfill\
\begin{tabular}{|l|l|}
\hline
coefficient&monomial\\
\hline
$1$&$\varphi_t-f_x\varphi_x-f\varphi_{xx}-h\varphi_x-k_x\xi-k_u\varphi$\\
$u_x$&$\xi_t+f_{xx}\xi+f_{xu}\varphi+f_x(\varphi_u-\xi_x)+2f_u\varphi_x+f(2\varphi_{xu}-\xi_{xx})+h(\varphi_u-\xi_x)+h_x\xi+h_u\varphi$\\
$u_t$&$\varphi_u-\eta_t+f_x\eta_x+f\eta_{xx}+h\eta_x$\\
$u_xu_t$&$\xi_u-f_x\eta_u-2f_u\eta_x-f\eta_{xu}-h\eta_u$\\
$u_t^2$&$\eta_u$\\
$u_x^2$&$-f_x\xi_u+f_{xu}\xi+f_{uu}\varphi+2f_u(\varphi_u-\xi_x)-h\xi_u+f(\varphi_{uu}-2\xi_{xu})$\\
$u_x^3$&$2f_u\xi_u+f\xi_{uu}$\\
$u_x^2u_t$&$2f_u\eta_u-f\eta_{uu}$\\
$u_{xx}$&$f_x\xi+f_u\varphi+f(\varphi_u-2\xi_x)$\\
$u_xu_{xt}$&$2f\eta_x$\\
$u_xu_{xx}$&$3f\xi_u$\\
$u_tu_{xx}$&$f\eta_u$\\
$u_{xx}u_x^2$&$2f\eta_u$\\
\hline
\end{tabular} \hfill\ (Table 1) \hfill\ \mbox{ }
\end{center}
By simplifying above equations we obtain:
\begin{eqnarray}\label{eq:5}
&&\eta=\eta(t),\qquad \xi=\xi(x,t),\qquad \varphi_u-\eta_t=0,\qquad f_x\xi+f_u\varphi+f(\varphi_u-2\xi_x)=0,\nonumber\\
&&\varphi_t-f_x\varphi_x-f\varphi_{xx}-h\varphi_x-k_x\xi-k_u\varphi=0,\qquad f_{xu}\xi+f_{uu}\varphi+2f_u(\varphi_u-\xi_x)=0,\\
&&\xi_t+f_{xx}\xi+f_{xu}\varphi+(f_x+h)(\varphi_u-\xi_x)+2f_u\varphi_x-f\xi_{xx}+h_x\xi+h_u\varphi=0.\nonumber
\end{eqnarray}
\section{Group classification in special cases}
In this section we consider four special case of modeling equation
and obtain differential invariants for them by using (\ref{eq:5})
and infinitesimal criterion method.
\paragraph*{A:} $f(x,u)=xu^{-1}$, $h(x,u)=-2/u$, $k(x,u)=au+b$, where $a,b$ are constant real
numbers. In this case we have:
$\eta=\frac{1}{a}e^{at}.e^{ac_1}+c_2$, $\xi= c_3\sqrt{x}$, and
$\varphi=u.e^{at}.e^{ac_1}$. As a result we find 3 independent
vector fields: $X_1=e^{at}\partial_t+ae^{at}\partial_u$,
$X_2=\partial_t$, $X_3= \sqrt{x}\partial_x$.
\paragraph*{B:} $f(x,u)=ax^4u$,
$h(x,u)={\frac{bx}{u}}$, $k(x,u)=xu$, where $a\neq0,b$ are real
numbers. In this case we have: $\eta=-c_1t+c_2$, $\xi=c_1x$,
$\varphi=-c_1u$. As a result we find 2 independent vector fields:
$X_1=-t\partial_t+x\partial_x-u\partial_u$, $X_2=\partial_t$.
\paragraph*{C:} $f(x,u)=ax\exp{(-u/b)}$, $h(x,u)=xu$,
$k(x,u)=c-bu$, where $a\neq,b\neq0,c$ are real constant numbers.
In this case we have: $\eta=c_1$, $\xi=-\frac{c_2}{c}x\exp{(bt)}$,
$\varphi=c_2\exp{(bt)}$. As a result we find 2 independent vector
fields:
$X_1=-\frac{1}{b}x\exp{(bt)}\partial_x+\exp{(bt)}\partial_u,$
$X_2=\partial_t$.
\paragraph*{D:} $f(x,u)=ax^2u$, $h(x,u)=xu$, $k(x,u)=u$, where $a$ is real nonzero
constant. In this case we have: $\eta=c_1$, $\xi=c_2x$,
$\varphi=0$. As a result we find 2 independent vector fields:
$X_1=\partial_t$, $X_2=x\partial_x$.

Similar to above, the reader can use above procedure for finding
her or him interested modeling equation where has form
(\ref{eq:11}), with interested $f,h$ and $k$.
\section{Resulted differential invariants}
In this section we obtain differential invariants for above
resulted symmetry groups in several major and complicated cases.
For example we compute differential equation for {\bf B}, $X_1$
and {\bf C}, $X_1$.
\paragraph*{B,} $X_1$: In this case we have following determination equation:
$ \frac{dx}{x}=\frac{dt}{-t}=\frac{du}{-u}$, and by solving this
equation we find: $xt=c_1$, $xu=c_2$, $u/t=c_3$; and we choose
$r=xt$ and $w=xu$ as independent invariants. (we note $u/x=w/r$
and as a result obtain from $r,w$.)
\paragraph*{C,} $X_1$: In this case we have following determination equation:
$ \frac{b dx}{x\exp{(bt)}}=\frac{dt}{0}=\frac{du}{\exp{(bt)}}$,
and by solving this equation we find : $t=c_1$, $c_2=u+b\ln{x}$;
and we choose $r=t$ and $w=u+b\ln{x}$ as independent invariants.
\begin{center}\small \mbox{ } \hfill\
\begin{tabular}{|l|l|l|}
 \hline
case& interested vector field&differential invariants\\
\hline
A&$X_1$& $r=x$, $w=u-at$\\
&$X_2$& $r=x$, $w=u$\\
&$X_3$& $r=t$, $w=u$\\
\hline
B&$X_1$&$r=xt$, $w=xu$\\
&$X_2$&$r=x$, $w=u$\\
\hline
C&$X_1$&$r=t$, $w=u+b\ln{x}$\\
&$X_2$&$r=x$, $w=u$\\
\hline
D&$X_1$&$r=x$, $w=u$\\
&$X_2$&$r=t$, $w=u$\\
\hline
\end{tabular} \hfill\ (Table 2) \hfill\ \mbox{ }
\end{center}
In the above table, $F$ is an arbitrary function.

In the sequel, we obtain reduced equation respect to specified
group symmetry with infinitesimal generator $X$, (solution of this
reduced equation called $X$-invariants solution of original
equation) by using resulted differential invariants in the above
table in two case.

For example, consider $u_t=(ax^4u)_x+bx/uu_x+xu$ (case B). By
considering $w=w(r)$, we find: $u_t=w_r$, $u_x=(xtw_r-w)/x^2$ and
$u_{xx}=(x^2(tw_r+xt^2w_{rr}-tw_r)-2x(xtw_r-w))/x^4$. By
substituting this values in the given equation, we find following
$X_1$-reduced equation:
\begin{eqnarray*}
w_r=b+(1-4a)w+(6a+ar^3)w^2+(4ar-2aw+awr-2arw)w_r+awr(r-a)w_{rr}
\end{eqnarray*}

As and second example, consider
$u_t=\big({\frac{ax}{\exp{(u/b)}}}\big)_x+xuu_x+c-bu$, By
considering  $w=w(r)$, we find: $u_t=w_r$, $u_x=-1/x$ and
$u_{xx}=1/x^2$. By substituting this values in the given equation,
we find following $X_1$-reduced equation:
\begin{eqnarray*}
&& w_r=c+ab-2a,
\end{eqnarray*}
%
%
\section{Some Applications}
The Kolomogorov-Petrovskii-Piskonov (KPP) equation, (See
\cite{[1],[7]})
\begin{eqnarray}
E(u)\equiv bu_t-u_{xx}+\gamma uu_x+f(u),\label{eq:6}
\end{eqnarray}
with ($b$,$\gamma$) real numbers, is encountered in
reaction-diffusion systems and prey-predator models. The optional
convection term $uu_x$ \cite{[1],[3]}) is quite important in
physical applications to prey-predator models.
\subsection{Classical symmetries and Differential invaiants}
If we let $b\neq0$, then we have following equation:
\begin{eqnarray}
u_t=\frac{1}{b}(u_{xx}-\gamma uu_x-f(u)),\label{eq:7}
\end{eqnarray}
By substituting this value in (\ref{eq:5}), we have following
results.
\paragraph*{Case I:} $b={\frac{\alpha\gamma}{\exp(\beta\alpha)}},
f(u)={\frac{(1/2)\gamma \kappa\alpha u}{\exp(\alpha\beta)}}+s$;
Where $\alpha,\beta,\kappa$ and s are arbitrary constants.

In this case we have:
\begin{eqnarray}
\xi=\frac{\exp(\alpha t)\exp(\alpha\beta)}{\alpha}+c_2,\qquad
\eta=c_1,\qquad \varphi=\kappa\exp(\alpha t),\label{eq:8}
\end{eqnarray}
For symmetry algebra we find:
\begin{eqnarray}
X_1=\partial_t \qquad X_2=\partial_x,\label{eq;9}
\end{eqnarray}
\paragraph*{Case II:} $b\neq{\frac{\alpha\gamma}{\exp(\beta\alpha)}},
f(u)\neq{\frac{(1/2)\gamma\kappa\alpha u}{\exp(\alpha\beta)}}+s$;
Where $\alpha,\beta,\kappa$ and s are arbitrary constants.

In this case we have:
\begin{eqnarray}
\xi=c_1\qquad  \eta=c_2,\qquad \varphi=0,\label{eq:8}
\end{eqnarray}
For symmetry algebra we find:
\begin{eqnarray}
X_1=\partial_t &&\qquad  X_2=\partial_x,\label{eq;9}
\end{eqnarray}
As a result we have following theorem:
\begin{thm}
Some exact solutions for modeling equation (\ref{eq:6}) invariant
under a translation group respect to $x$ and some solutions of
this equation invariant under translation respect to $t$.
\end{thm}
\subsection{Similarity solutions}
In this subsection we find similarity solution of equation
(\ref{eq:7}) by using above resulted symmetry algebra.
\paragraph*{similarity solution respect to $X=\partial_t$.}
In this case we have following equation as $X$-reduced equation:
\begin{eqnarray}
w_{rr}-\gamma ww_r-f(w)=0,\label{eq:12}
\end{eqnarray}
If we solve equation (\ref{eq:12}) with MAPLE, then we find:
$w(x)=c$. Where
\begin{eqnarray}
\frac{d}{dc}{F(c)}F(c)-\gamma(c F(c))-f(c)=0\hspace{.5cm}
\mbox{or} \hspace{0.5cm} \frac{d}{dr}w(r)=F(c), \hspace{.5cm}
\mbox{or}\hspace{.5cm} r=\int \frac{1}{F(c)}dc+C
\end{eqnarray}
Where $F$ is arbitrary function with specified arguments and $c,C$
are arbitrary constants.
\paragraph*{similarity solution respect to $Y=\partial_x$.}
In this case we have following equation as $Y$-reduced equation:
\begin{eqnarray}
w_{r}+\frac{1}{b}f(w)=0,\label{eq:13}
\end{eqnarray}
If we solve equation (\ref{eq:13}) with MAPLE, then we find
following solution:
\begin{eqnarray}
x-\int^{w(x)}\frac{b}{f(c_1)}dc_1+c_2=0,
\end{eqnarray}
Where $c_1$ and $c_2$ are arbitrary constants.
%
\section*{Conclusion}
In this paper first we find system of equations to finding
symmetry group and symmetry algebra for (G-RDC) equation, then
obtain these symmetry groups in several special cases and at the
end we establish symmetry classification for KPP equation by using
group classification of (G-RDC) equation and we find its
similarity solution respect to resulted symmetry algebra.
%

%
\end{document}